\documentclass[10pt,twoside]{article}

\usepackage[english]{babel}
\usepackage[utf8]{inputenc}

\usepackage{amsthm, amsmath, amssymb, bbm}
\RequirePackage[all,cmtip]{xy}

\usepackage{tikz}
\usetikzlibrary{calc}

\usepackage{graphicx}

\usepackage[normalem]{ulem}
\newcommand\redsout{\bgroup\markoverwith{\textcolor{red}{\rule[0.5ex]{2pt}{1pt}}}\ULon}
\newcommand\redssout{\bgroup\markoverwith{\textcolor{red}{\rule[0ex]{2pt}{.5pt}}}\ULon}

\makeatletter
\newcommand{\email}[1]{{\footnotesize {\it E-mail address}: {\tt #1\vphantom{y}}}}

\def\address#1{\begingroup\setlength{\parindent}{0pt}%
\vspace*{15pt}
\fontsize{10}{8}\selectfont\sc #1
\endgroup}

\@addtoreset{equation}{section}

\theoremstyle{plain}

\newtheorem{theorem}{Theorem}[section]
\newtheorem{proposition}[theorem]{Proposition}

\theoremstyle{definition}

\newtheorem{example}[theorem]{Example}

\newtheorem{remark}[theorem]{Remark}
\newtheorem{problem}[theorem]{Problem}

\renewcommand{\emptyset}{\varnothing}

\newcommand{\C}{\mathbb{C}}

\newcommand{\R}{\mathbb{R}}

\newcommand{\Z}{\mathbb{Z}}

\newcommand{\op}{\mathrm{op}}

\DeclareMathOperator{\id}{id}

\newcommand{\BDC}{\mathbf{D}^{\mathrm{b}}}

\newcommand{\DSum}{\bigoplus}
\newcommand{\dsum}[1][]{\mathbin{\oplus_{#1}}}

\renewcommand{\to}[1][]{\xrightarrow[#1]{}}

\newcommand{\isoto}[1][]{\xrightarrow[#1]{\sim}}

\newcommand{\Endo}[1][]{\mathrm{End}_{\raise1.5ex\hbox to.1em{}#1}}

\newcommand{\Hom}[1][]{\mathrm{Hom}_{\raise1.5ex\hbox to.1em{}#1}}

\newcommand{\RHom}[1][]{\mathrm{RHom}_{\raise1.5ex\hbox to.1em{}#1}}

\newcommand{\Ext}[2][]{\mathrm{Ext}_{\raise1.5ex\hbox to.1em{}#1}^{#2}}

\newcommand{\THom}[1][]{\mathrm{THom}_{\raise1.5ex\hbox to.1em{}#1}}

\newcommand{\Mod}{\mathrm{Mod}}

\newcommand{\Tens}[1][]{\mathbin{\otimes_{\raise1.5ex\hbox to-.1em{}#1}}}

\newcommand{\LTens}[1][]{\mathbin{\otimes_{\raise1.5ex\hbox to-.1em{}#1}^{L}}}

\newcommand{\Tor}[2][]{\mathrm{Tor}^{\raise1.5ex\hbox to.1em{}#1}_{#2}}

\def\she{\mathcal{E}}

\def\shl{\mathcal{L}}
\def\shm{\mathcal{M}}

\newcommand{\sect}{\Gamma}

\renewcommand{\hom}[1][]{{\mathcal{H}om}_{\raise1.5ex\hbox to.1em{}#1}}

\newcommand{\rhom}[1][]{{R\mathcal{H}om}_{\raise1.5ex\hbox to.1em{}#1}}

\newcommand{\ext}[2][]{{\mathcal{E}xt}_{\raise1.5ex\hbox to.1em{}#1}^{#2}}

\newcommand{\thom}[1][]{{T\mathcal{H}om}_{\raise1.5ex\hbox to.1em{}#1}}

\newcommand{\tens}[1][]{\mathbin{\otimes_{\raise1.5ex\hbox to-.1em{}#1}}}

\newcommand{\ltens}[1][]{\mathbin{\otimes_{\raise1.5ex\hbox to-.1em{}#1}^{L}}}

\newcommand{\tor}[2][]{{\mathcal{T}or}^{\raise1.5ex\hbox to.1em{}#1}_{#2}}

\newcommand\etens{\mathbin{\boxtimes}}

\newcommand{\reim}[1]{{R#1}_!}

\newcommand{\opb}[1]{#1^{-1}}

\newcommand{\epb}[1]{#1^{!}}

\newcommand{\GHom}[1][]{\mathrm{GHom}_{\raise1.5ex\hbox to.1em{}#1}}

\newcommand{\GExt}[2][]{\mathrm{GExt}_{\raise1.5ex\hbox to.1em{}#1}^{#2}}

\newcommand{\FHom}[1][]{\mathrm{FHom}_{\raise1.5ex\hbox to.1em{}#1}}

\newcommand{\ghom}[1][]{{\mathcal{GH}om}_{\raise1.5ex\hbox to.1em{}#1}}

\newcommand{\gext}[2][]{{\mathcal{GE}xt}_{\raise1.5ex\hbox to.1em{}#1}^{#2}}

\newcommand{\fhom}[1][]{{\mathcal{FH}om}_{\raise1.5ex\hbox to.1em{}#1}}

\newcommand{\D}{\mathcal{D}}

\renewcommand{\O}{\mathcal{O}}

\newcommand{\Db}{\mathcal{D}b}

\newcommand{\detens}{\mathbin{\underline{\boxtimes}}}

\def\absdoim#1{\underline{#1}_*}
\def\reldoim[#1]#2{\underline{#2}_{|{#1}*}}
\def\doim{\@ifnextchar [{\reldoim}{\absdoim}}

\def\absdeim#1{\underline{#1}_*}
\def\reldeim[#1]#2{\underline{#2}_{|{#1}*}}
\def\deim{\@ifnextchar [{\reldeim}{\absdeim}}

\def\absdopb#1{\underline{#1}^{-1}}
\def\reldopb[#1]#2{\underline{#2}_{|{#1}}^{-1}}
\def\dopb{\@ifnextchar [{\reldopb}{\absdopb}}

\newcommand{\hol}{\mathrm{hol}}

\newcommand{\reghol}{\mathrm{r-hol}}

\newcommand{\defeq}{\mathbin{:=}}
\newcommand{\eqdef}{\mathbin{=:}}

\renewcommand{\Re}{\operatorname{Re}}

\newcommand{\Cc}{{\C\text-\mathrm{c}}}
\newcommand{\Rc}{{\R\text-\mathrm{c}}}
\newcommand{\sol}[1][]{\mathcal Sol_{#1}}

\renewcommand{\BDC}{\mathrm{D}^{\mathrm{b}}}

\renewcommand{\reim}[1]{{\mathrm{R}#1}_!}

\newcommand{\sub}{{\operatorname{sub}}}
\newcommand{\temp}{{\operatorname{temp}}}
\newcommand{\reeim}[1]{{\mathrm{R}#1}_{!!}}

\newcommand{\enh}{{\operatorname{enh}}}
\newcommand{\enhsub}{{\operatorname{enh-sub}}}
\newcommand{\Temp}{{\operatorname{T}}}
\newcommand{\ctens}{\mathbin{\overset+\otimes}}
\newcommand{\indlim}[1][]{\mathop{\text{\rm``$\operatorname{colim}$''}}\limits_{#1}}

\newcommand{\dtens}[1][]{\mathbin{\otimes_{\raise1.5ex\hbox to-.1em{}#1}^{\mathsf{D}}}}
\renewcommand{\detens}[1][]%
{\mathbin{\boxtimes_{\raise1.5ex\hbox to-.1em{}#1}^{\mspace{2mu}\mathsf{D}}}}
\renewcommand{\doim}[1]{{\mathsf{D}#1}_*\mspace{1mu}}
\renewcommand{\dopb}[1]{{\mathsf{D}#1}^{\mspace{1mu}*}}

\renewcommand{\reghol}{\text{\rm reg-hol}}
\renewcommand{\Cc}{{\C\text{\rm -cons}}}
\renewcommand{\Rc}{{\R\text{\rm -cons}}}
\makeatother

\begin{document}

\title{On the irregular Riemann-Hilbert correspondence}

\author{by Andrea D'Agnolo and Masaki Kashiwara}

\date{}

\maketitle

\markboth{\hfill{\rm Andrea D'Agnolo, Masaki Kashiwara} \hfill}{\hfill {\rm On the irregular Riemann-Hilbert correspondence \hfill}}

\section*{Introduction}

The original Riemann-Hilbert problem asks to find a Fuchsian ordinary differential equation with prescribed singularities and monodromy in the complex line. In the early 1980's Kashiwara \cite{Kas84} solved a generalized version of the problem, valid on complex manifolds of any dimension. He presented it as a correspondence between regular holonomic $\D$-modules and perverse sheaves.

The analogous problem where one drops the regularity condition remained open for about thirty years.  We solved it in the paper \cite{DK16} that just received a 2024 Frontiers of Science Award.  Our construction requires in particular an enhancement of the category of perverse sheaves.

Here\footnote{This is a written account of a talk given by the first named author at the International Congress of Basic Sciences on July 2024 in Beijing.}, using some examples in dimension one, we wish to convey the gist of the main ingredients used in our work.
For another, more detailed, presentation we refer to \cite{Kas16}.

\section{Original formulation}\label{sec:ODE}

The Riemann-Hilbert problem is the 21st in the list that
Hilbert presented at the Paris conference of the International Congress of Mathematicians in 1900. It was probably motivated by the work \cite{Rie57} of
Riemann in 1857 on the hypergeometric functions. The problem is stated as follows:

\begin{problem}\label{pro:RH}
To find a Fuchsian ordinary differential equation with prescribed singularities and monodromy.
\end{problem}

Here is a toy model. For $\lambda\in\C$, the equation $(z\frac d{dz}-\lambda)u=0$ has $u=z^\lambda$ as a basis of solutions. As $u$ has monodromy $e^{2\pi i \lambda}$ around the origin, we see that Problem~\ref{pro:RH} is of a transcendental nature. There are no other Fuchsian operators with the same monodromy, up to meromorphic gauge equivalence.

For non Fuchsian equations, like the Airy equation, Stokes \cite{Sto57} observed in the same year 1857 a phenomenon beyond monodromy. This phenomenon now bears his name, and we recall it in Section~\ref{se:stokes}.
The case where one drops the Fuchsian condition is addressed by the Riemann-Hilbert-Birkhoff problem:

\begin{problem}\label{pro:RHB}
To find an ordinary differential equation with prescribed singularities, monodromy and Stokes data.
\end{problem}

\section{Modern incarnation}

Let $(X,\O_X)$ be a complex manifold, and denote by $\D_X$ the sheaf of algebras of linear partial differential operators of finite order.
In 1984 Kashiwara \cite{Kas84} solved a wide generalization of Problem~\ref{pro:RH} with the following theorem.

\begin{theorem}\label{th:Kas84}
There is an equivalence of abelian categories, compatible with Grothendieck's operations,
\[
\xymatrix@C=4em@R=0em{
\left\{
\text{regular holonomic $\D_X$-modules}
\right\}
\ar@<.8ex>[r]^-{\sol}_-\sim
&
\left\{
\text{perverse sheaves on $X$}
\right\}^\op
\ar@<.9ex>[l]
\\
\shm \ar@{|->}[r]
& \rhom[\D_X](\shm,\O_X).
}
\]
\end{theorem}

Holonomic $\D_X$-modules, formerly known as maximally overdetermined systems, are a higher dimensional analogue of ordinary differential equations. Similarly, regularity is a higher dimensional analogue to the Fuchsian condition. Perverse sheaves, despite their name, are not sheaves but objects of the derived category. It thus makes sense to compute the derived hom functor, used in the definition of $\sol$. This functor is a higher dimensional analogue of computing holomorphic solutions for an ordinary differential equation. As indicated by the arrow pointing to the left, an explicit quasi-inverse to $\sol$ was also given in \cite{Kas84}.

Even without the regularity condition, Kashiwara \cite{Kas75} proved back in 1975 that the solution functor still takes values in perverse sheaves.
Then, the first problem to establish an equivalence for possibly irregular holonomic $\D$-modules was
to come up with an enhancement of that target category. Such a problem remained unsettled for about thirty years until 2013, when we submitted the paper \cite{DK16} that just received a 2024 Frontiers of Science Award. Here is our solution:

\begin{theorem}\label{th:RHirr}
There is an equivalence\footnote{\label{foot:disc}Our statement was not so clear cut, as the category of perverse enhanced subanalytic sheaves still lacks a purely topological description. We will comment on this in Section~\ref{sec:RH}.} of abelian categories, compatible with Grothendieck's operations,
\[
\xymatrix@C=4em@R=0em{
\left\{
\text{holonomic $\D_X$-modules}
\right\}
\ar@<.8ex>[r]^-{\sol^\Temp}_-\sim
&
\left\{
\begin{minipage}{12.5em}
perverse enhanced subanalytic\newline
sheaves on $X$
\end{minipage}
\right\}^\op .
\ar@<.9ex>[l]
}
\]
\end{theorem}

Some explanations are in order, both for the adjectives ``enhanced'' and ``subanalytic'', as well as for the tempered enhancement $\sol^\Temp$ of the solution functor. We will do it through some significant examples in dimension one.
For this, let us unwind back to the setting of Section~\ref{sec:ODE}.

\section{The case of meromorphic connections}\label{Deligne}

Denote by $\O_X(*D)$ the sheaf of meromorphic functions with possible poles on a hypersurface $D\subset X$.
In 1970 Deligne \cite{Del70} proved the following statement, which does not use derived categories.

\begin{theorem}\label{th:Del}
There is an equivalence of abelian categories
\[
\xymatrix@C=4em{
\left\{
\text{regular flat $\O_X(*D)$-connections}
\right\}
\ar[r]^-{\sol|_{X\setminus D}}_-\sim
&
\left\{
\text{local systems on $X\setminus D$}
\right\}^\op
}.
\]
\end{theorem}

In the terminology of Section~\ref{sec:ODE}, $D$ are the prescribed singularities, and the monodromy is encoded as the action of the fundamental group of $X\setminus D$ on the stalks of a given local system.

Note that Theorem~\ref{th:Del} is encompassed by  Theorem~\ref{th:Kas84}. In fact, the category of (regular) flat $\O_X(*D)$-connections is but the full subcategory of (regular) holonomic $\D_X$-modules $\shm$ such that:
\begin{equation}
\label{eq:flatconn}
\operatorname{SingSupp}(\shm)\subset D, \quad \shm\isoto\shm\tens[\O]\O_X(*D).
\end{equation}
This amounts to ask that $\shm|_{X\setminus D}$ is a (regular) flat $\O_{X\setminus D}$-connection, and that $\shm$ is stable by localization at $D$.

For $f\in\O_X(*D)$, a basic example of meromorphic connection is
\begin{equation}\label{eq:E}
\she_{X\setminus D}^f\defeq (\O_X(*D),\, d+df),
\end{equation}
which is regular if and only if $f\in\O_X$.

\section{The case of ODEs}

Let $X\subset\C_z$ be an open neighborhood of the origin $z=0$, and set $D= \{0\}$.
Consider a differential operator
\[
P=a_m(z)\left(\tfrac d{dz}\right)^m+\cdots+a_1(z)\tfrac d{dz}+a_0(z),
\qquad a_j\in\O_X,
\]
and assume that its top degree coefficient $a_m$ vanishes only at $z=0$.
Then the associated $\D_X$-module
$\shm \defeq (\D_X/\D_X P)\tens[\O]\O_X(*D)$
satisfies conditions \eqref{eq:flatconn}.
In this case,
$\sol(\shm)|_{X\setminus D} = \{ u\in \O_{X\setminus D}\colon Pu=0\}$
is the rank $m$ local system of usual holomorphic solutions.

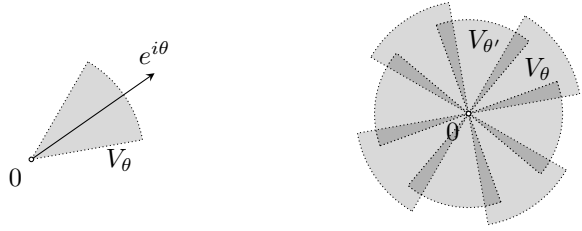
\begin{figure}
\hfill
\begin{tikzpicture}[scale=.5,baseline={(current bounding box.center)}]
\draw[color=black,fill=black!15,densely dotted] (60:3) -- (0,0) -- (10:3)  node[below left]{$V_\theta$}
arc[start angle=10, end angle=60, radius=3]
  -- cycle ;
\draw[-stealth](0,0) -- (35:4) node[above]{$e^{i\theta}$};
\filldraw[color=black,fill=white] (0,0) circle (.06) node[below left]{$0$} ;
\end{tikzpicture}
\hfill
\begin{tikzpicture}[scale=.5,baseline={(current bounding box.center)}]
\foreach \x in {0,...,3}
	{\draw[color=black,fill=black!15,densely dotted] (60+\x*90:3) -- (0,0) -- (10+\x*90:3) arc[start angle=10+\x*90, end angle=60+\x*90, radius=3]
  -- cycle ;}
\foreach \x in {0,...,3}
	{\draw[color=black,fill=black!15,densely dotted] (20+\x*90:2.5) -- (0,0) -- (-40+\x*90:2.5) arc[start angle=-40+\x*90, end angle=20+\x*90, radius=2.5]
  -- cycle ;}
\foreach \x in {0,...,3}
	{\draw[color=black,fill=black!30,densely dotted] (20+\x*90:2.5) -- (0,0) -- (10+\x*90:2.5) arc[start angle=10+\x*90, end angle=20+\x*90, radius=2.5]
  -- cycle ;}
\foreach \x in {0,...,3}
	{\draw[color=black,fill=black!30,densely dotted] (60+\x*90:2.5) -- (0,0) -- (50+\x*90:2.5) arc[start angle=50+\x*90, end angle=60+\x*90, radius=2.5]
  -- cycle ;}
\draw (35:3)  node[below left]{$V_\theta$} ;
\draw (80:2.5)  node[below]{$V_{\theta'}$} ;
\filldraw[color=black,fill=white] (0,0) circle (.06) node[below left]{$0$} ;
\end{tikzpicture}
\hfill
\phantom{a}
\caption{sectorial neighborhoods of directions emanating from $z=0$.}\label{fig:sectors}
\end{figure}

When $P$ is not necessarily Fuchsian, additional information besides monodromy can be obtained by considering formal or asymptotic solutions to $Pu=0$. They are described by this classical result due to Hukuhara, Levelt and Turrittin:

\begin{proposition}\label{pro:HLT}
$Pu=0$ has bases of
\begin{itemize}
\item[(i)] formal solutions $\widehat u_1,\dots,\widehat u_m$ at $z=0$, of the form
\[
\widehat u_j = z^{\lambda_j} e^{f_j(z)} \sum_{k=0}^{m-1}\widehat a_{jk}(z)(\log z)^k,
\]
where $\lambda_j\in\C$, $f_j\in z^{-1/d}\,\C[z^{-1/d}]$ for a ramification index $d\in\Z_{\geq 1}$,
and $\widehat a_{jk}(z)
= \sum_{\ell= 0}^{+\infty}a_{jk\ell} \,z^{\ell/d}$ is a formal power series in $z^{1/d}$.
\item[(ii)] asymptotic solutions $u_1^\theta,\dots,u_m^\theta\in\O_X(V_\theta)$,
for any $\theta\in\R$, with $u_j^\theta\underset{(*)}\sim \widehat u_j$. Here, $V_\theta$ is a small enough open sectorial neighborhood of the direction $e^{i\theta}$, as in the left of Figure~\ref{fig:sectors}, and the relation $(*)$ means
\begin{align*}
{}^\forall N\in\Z_{>0}\ {}^\exists C>0\ {}^\forall z\in V_\theta\colon
\quad
|u_j(z) - \widehat u_j^N(z)|
&\leq C\,|e^{f_j(z)} z^{\lambda_j+N/d}| \\
&= C\,e^{\Re f_j(z)}|z^{\lambda_j+N/d}|.
\end{align*}
In these estimates, $\widehat u_j^N$ denotes the same expression as $\widehat u_j$, but with the formal series $\widehat a_{jk}$ replaced by their partial sum $\widehat a^N_{jk}(z)
= \sum_{\ell= 0}^N a_{jk\ell} \,z^{\ell/d}$. Moreover, in order for $\widehat u_j^N$ to make sense as a function, one has to fix on $V_\theta$ a determination of $\log z$ (and hence of $z^{\lambda_j}$ and $z^{1/d}$).
\end{itemize}
\end{proposition}

The above asymptotic basis of solutions $u_j^\theta\sim \widehat u_j$ is not unique. In fact,
\begin{equation}\label{eq:ujtheta}
\Re f_k<\Re f_j \text{ on } V_\theta
\implies u_j^\theta+s\, u_k^\theta \sim \widehat u_j
\quad  ^\forall s\in\C,
\end{equation}
due to the rapid decay of $e^{f_k(z)-f_j(z)}$ at $z=0$ along $V_\theta$.

\section{Stokes data}\label{se:stokes}

Let us keep notations as in the previous section.
For $j=1,\dots,m$, the exponents $f_j$ in the expression of $\widehat u^j$ are called \emph{exponential factors} of $P$ at $z=0$. They are the first ingredient of the Stokes data mentioned in Problem~\ref{pro:RHB}. The second and last ingredient are \emph{Stokes matrices}. They are obtained by the following recipe:
\begin{enumerate}
\item fix a basis of formal solutions $\widehat u_j$ of $Pu=0$ at $z=0$,
\item surround $z=0$ with finitely many sectors $V_\theta$, as on the right of Figure~\ref{fig:sectors},
\item on each such sector $V_\theta$, choose a basis of asymptotic solutions $u_j^\theta\sim \widehat u_j$.
\end{enumerate}
To each overlapping pair of sectors $V_\theta\cap V_{\theta'}\neq\emptyset$ is attached a
Stokes matrix. This is the matrix $S_{\theta'\theta}=(s_{jk})$ of base change for asymptotic solutions, given by:
\[
u^{\theta'}_j=\sum_{k=1}^m s_{jk}\, u^\theta_k \qquad\text{on } V_\theta\cap V_{\theta'}.
\]
Note that \eqref{eq:ujtheta} gives
\begin{equation}\label{eq:Stokes}
s_{jk} \neq 0 \implies \Re f_k \leq \Re f_j \text{ on }V_\theta\cap V_{\theta'}.
\end{equation}

As monodromy is encoded in local systems, so Stokes data can be encoded in ``Stokes local systems''\footnote{As the quotes indicate, this notion of ``Stokes local system'' has some drawbacks. We will upgrade it to an unquoted version in Section~\ref{sec:SLS}.}. Let us present this notion in the case of the Airy equation considered by Stokes~\cite{Sto57}.

\begin{example}\label{ex:Airy}
Let $X=\mathbb{P}^1=\C_z\cup\{\infty\}$, and set $D=\{\infty\}$.
The Airy operator $P=\left(\frac d{dz}\right)^2-z$ has no singularity at finite distance, and has a non Fuchsian singularity at $z=\infty$.
There, its exponential factors are $f_{1,2}(z) = \pm \frac23 z^{3/2}$, with ramification index $2$.
In order to keep track of their growth, Stokes~\cite[page 116]{Sto57} drew the picture on the left of Figure~\ref{fig:Stokes}.
As he explains, it is a polar plot of the multivalued functions $\rho+\sigma\Re f_{1,2}(z)$, for some $\rho,\sigma\in\R_{>0}$, with $|z|=R\gg0$. (Thus, the dotted circle has radius $\rho$.)

\begin{figure}
\hfill
\begin{tikzpicture}[scale=.7,baseline=(O.base)]
\draw[densely dotted] (0,0)  coordinate (O) circle (2) ;
\def\rho{(2+cos(3*(\t)/2)/2.9)}
\begin{scope}
\pgfplothandlerlineto
\pgfplotfunction{\t}{0,...,720}
     {\pgfpointxy {cos(\t)*\rho}{sin(\t)*\rho}}
     \pgfusepath{stroke}
\end{scope}
\filldraw (0:1.66) circle (.03) node[left] {$\scriptstyle a$} ;
\filldraw (0:2.35) circle (.03) node[right] {$\scriptstyle B$} ;
\filldraw (120:1.66) circle (.03) node[below right] {$\scriptstyle c$} ;
\filldraw (120:2.35) circle (.03) node[above left] {$\scriptstyle A$} ;
\filldraw (240:1.66) circle (.03) node[above right] {$\scriptstyle b$} ;
\filldraw (240:2.35) circle (.03) node[below left] {$\scriptstyle C$} ;
\end{tikzpicture}
\hfill
\begin{tikzpicture}[scale=1.3,baseline=(O.base)]
\begin{scope}
    \clip(0,-.5) rectangle (4.5,.8);
	\draw[color=black,fill=black!15,dashed] (150/80,-1) rectangle (170/80,1) ;
	\draw[color=black,fill=black!15,dashed] (230/80,-1) rectangle (250/80,1) ;
	\draw[color=black,fill=black!15,dashed] (300/80,-1) rectangle (320/80,1) ;
\end{scope}
\draw (160/80,-.5)  node[below]{\tiny$\begin{pmatrix}*\ 0\\ *\ *\end{pmatrix}$\rlap,};
\draw (240/80,-.5)  node[below]{\tiny$\begin{pmatrix}*\ *\\ 0\ *\end{pmatrix}$\rlap,};
\draw (310/80,-.5)  node[below]{\tiny$\begin{pmatrix}*\ 0\\ 0\ *\end{pmatrix}$};
\draw (.7,-.55)  node[below]{$\scriptstyle{\operatorname{End}\left(E_M^{\varphi_1}\dsum E_M^{\varphi_2}\right)\ =}$};
\draw[densely dotted,-stealth] (0,0)  coordinate (O) -- (4.5,0)  node[right]{$\scriptstyle\theta$};
\draw[densely dotted,-stealth] (10/80,-.5) -- (10/80,.8)  node[below right]{$\scriptstyle t$};
\draw (4.5,1/3.5) node[above]{$\scriptstyle{-\varphi_1}$};
\draw (4.5,-1/3.5) node[below]{$\scriptstyle{-\varphi_2}$};
\def\rho{(cos(3*(\t-10)/2)/3.5)}
\begin{scope}
	\pgfplothandlerlineto
	\pgfplotfunction{\t}{0,...,360}
     {\pgfpointxy {\t/80}{\rho}}
     \pgfusepath{stroke}
\end{scope}
\begin{scope}
	\pgfplothandlerlineto
	\pgfplotfunction{\t}{0,...,360}
     {\pgfpointxy {\t/80}{-\rho}}
     \pgfusepath{stroke}
\end{scope}
\end{tikzpicture}
\hfill
\phantom{a}
\caption{plots of the real part of Airy's exponential factors on $|z|=R\gg0$.}\label{fig:Stokes}
\end{figure}
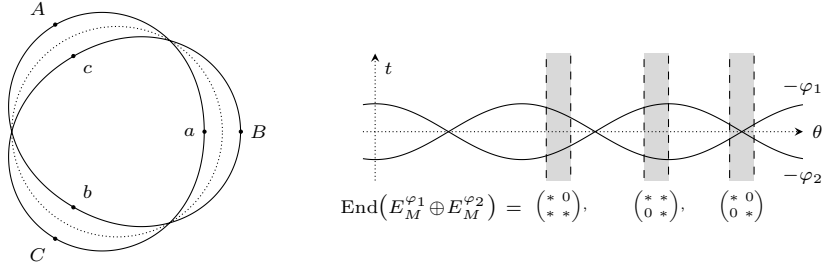

As $f_{1,2}$ are multivalued, let us consider the universal cover of $|z|=R$ given by $q\colon M=\R_\theta\to X$, $\theta\mapsto R\,e^{i\theta}$.
Then, $\varphi_{1,2}=q^*\Re f_{1,2}$ have a natural determination given by $\varphi_{1,2}(\theta)=\pm \frac23 R^{3/2}\cos(3\theta/2)$.
Consider the \emph{exponential} sheaf on $M\times\R_t$
\begin{equation}\label{eq:EM}
E_M^{\varphi_j} \defeq \C_{\{t+\varphi_j(\theta)\geq0\}},
\end{equation}
i.e.\ the extension by zero of the constant sheaf on the closed epigraph of $-\varphi_j$.

A ``Stokes local system'' on $M$ with exponential factors $\varphi_{1,2}$ is a sheaf $K$ on $M\times\R$ which, locally on $M$, is isomorphic to $E_M^{\varphi_1}\dsum E_M^{\varphi_2}$.
In other words, there is a cover $\{W_\theta\}_{\theta\in I}$ of $M$ such that $K|_{W_\theta\times\R}\simeq (E_M^{\varphi_1}\dsum E_M^{\varphi_2})|_{W_\theta\times\R}$ for any $\theta\in I$.
Thus, $K$ can be reconstructed from the induced automorphisms of $E_M^{\varphi_1}\dsum E_M^{\varphi_2}$ on overlapping vertical bands $(W_\theta\cap W_{\theta'})\times\R$.
As indicated on the right of Figure~\ref{fig:Stokes}, where the vertical bands are grayed out, such automorphisms satisfy the same condition \eqref{eq:Stokes} as Stokes matrices. In fact, assuming $W_\theta\cap W_{\theta'}$ connected, if there is a non zero morphism $u\colon E_M^{\varphi_j} \to E_M^{\varphi_k}$, then the epigraph of $-\varphi_k$ must be included in that of $-\varphi_j$, and hence $\varphi_k\leq\varphi_j$. In this case, $u=s_{jk}\,u_1$ for some $s_{jk}\in\C$, where $u_1\colon E_M^{\varphi_j} \to E_M^{\varphi_k}$ is the canonical morphism induced by the inclusion.
\end{example}

This motivates the introduction of enhanced sheaves, as a natural framework for this kind of constructions.

\section{Enhanced sheaves}

The theory of (what we call) enhanced sheaves was introduced by Tamarkin \cite{Tam18} in 2008\footnote{The paper is from 2018, well after ours, but the preprint was posted on the ArXiv in 2008.}, with motivations coming from symplectic topology.

\medskip
Let $M$ be a real analytic manifold, and consider the maps
\[
\xymatrix@C=2em{M\times\R^2 \ar[rr]^-{p_1,p_2,\mu} && M\times\R \ar[r]^-\pi & M},
\]
where $p_1,p_2,\pi$ denote the natural projections, and $\mu(x,t_1,t_2)=(x,t_1+t_2)$.
On the bounded derived category $\BDC(\C_{M\times\R})$ of sheaves on $M\times\R_t$, consider the additive convolution
\[
K_1\ctens K_2 \defeq \reim\mu(\opb p_1 K_1\tens\opb p_2 K_2).
\]
Note that $\C_{\{t=0\}}$ is a unit for $\ctens$.

Let $U\subset M$ be an open subanalytic subset, and $\varphi\colon U \to \R$ a continuous subanalytic\footnote{It means that the graph of $\varphi$  is subanalytic in $M\times\R$.} function. Generalizing \eqref{eq:EM}, consider the exponential sheaf on $M\times\R_t$
\begin{equation}\label{eq:EU}
E_U^\varphi\defeq\C_{\{t+\varphi(x)\geq 0\}},
\end{equation}
where now the epigraph of $-\varphi$, given by $\{(x,t)\colon x\in U,\ t+\varphi(x)\geq0\}$, is only  locally closed in $M\times\R$.
Note that $E_{U_1}^{\varphi_1}\ctens E_{U_2}^{\varphi_2}\simeq E_{U_1\cap U_2}^{\varphi_1+\varphi_2}$. In particular, $E^0_M=\C_{\{t\geq 0\}}$ is an idempotent for $\ctens$.

The derived category of enhanced sheaves on $M$ is defined by
\[
\BDC(\C^\enh_M) \defeq \{K\in\BDC(\C_{M\times\R})\colon K\ctens E_M^0\isoto K\},
\]
where the isomorphism is induced by the natural morphism $\C_{\{t\geq0\}}\to\C_{\{t=0\}}$.
It is a tensor category with respect to $\ctens$, with $E_M^0$ as unit.

External Grothendieck's operations for enhanced sheaves, relative to a morphism $r\colon M\to N$, are the corresponding operations for usual sheaves, relative to $\tilde r\defeq r\times\id_\R$. For example, the proper direct image of $K$ is $\reim{\tilde r}K$.

There is a fully faithful embedding, compatible with operations:
\begin{equation}\label{eq:suben}
\epsilon\colon\BDC(\C_M)\hookrightarrow \BDC(\C^\enh_M), \qquad F\mapsto\opb\pi F\ctens E_M^0.
\end{equation}

Enhanced sheaves are the correct framework for ``Stokes local systems''. However, this notion is not satisfactory in our quest for an irregular Riemann-Hilbert correspondence:

\medskip\noindent \textbf{Caveat I.}\quad
Let $X=\C_z$ and $D=\emptyset$. For $0\neq a\in\O_X$, let $P=\frac d{dz}-a'(z)$ be the operator associated with the regular connection $\she_X^a$ from \eqref{eq:E}.
Since both monodromy and Stokes data are trivial for $P$, its corresponding ``Stokes local system'' is $E_X^{\Re a}$.
However, one has $\she_X^a\simeq\she_X^0=\O_X$ and $E_X^{\Re a}\not\simeq E_X^0=\epsilon\,\C_X$.

\section{Subanalytic sheaves}

Subanalytic sheaves were introduced by Kashiwara-Schapira~\cite{KS01} in 2001, as a special case of their theory of ind-sheaves. We use them in our work to encode temperedness conditions for holomorphic functions.

\medskip
Let $M$ be a real analytic manifold.
Here is a definition of subanalytic sheaves, parallel to that of usual sheaves:
\begin{itemize}
\item Let $\operatorname{Op}_M$ be the category of open subsets of $M$.
The category $\Mod(\C_M)$ of usual sheaves on $M$ is the full subcategory of presheaves on $\operatorname{Op}_M$ satisfying the usual gluing conditions on covers.
\item Let $\operatorname{Op}_M^\sub\subset \operatorname{Op}_M$ be the full subcategory of  \emph{relatively compact, subanalytic} subsets.
The category $\Mod(\C_M^\sub)$ of subanalytic sheaves on $M$ is the full subcategory of presheaves on $\operatorname{Op}_M^\sub$ satisfying the usual gluing conditions on \emph{finite} covers.
\end{itemize}

There are natural functors\footnote{As in \eqref{eq:suben}, we denote by $\hookrightarrow$ a functor which is fully faithful.}
\[
\xymatrix@C=4em{
\Mod(\C_M)
\ar@<.5ex>@{^(->}[r]^-\iota
\ar@<-1.5ex>[r]_-\beta
& \Mod(\C_M^\sub) \ar@<.5ex>[l]|-\alpha ,
}
\]
defined as follows: $\iota$ restricts a usual presheaf to $\operatorname{Op}_M^\sub$;
$\alpha$ is the usual sheafification\footnote{This is well defined since the open subsets in $\operatorname{Op}_M^\sub$ form a basis for the topology of $M$} applied to a subanalytic presheaf;
and $\beta(F)$ is the subanalytic sheafification of $\operatorname{Op}_M^\sub \owns U\mapsto F(\overline U)$.

Subanalytic sheaves have analogous operations as those for usual sheaves, but not all of them commute with $\iota$. In particular, $\operatorname{colim}$ and $r_!$ do not commute with $\iota$.
Following \cite{KS01}, we denote by $\mathop{\text{\rm``$\operatorname{colim}$''}}$ and $r_{!!}$ the corresponding operations for subanalytic sheaves.
One has
\begin{equation}\label{eq:ilim}
\sect(U;\indlim[j]\iota\,F_j)\simeq\mathop{\operatorname{colim}}\limits_j\sect(U;F_j),
\qquad\text{for } U\in\operatorname{Op}_M^\sub.
\end{equation}

The subanalytic sheaf of \emph{tempered distributions} from \cite{KS01} is given by
\[
\Db^\temp_M(U) \defeq \operatorname{image}\bigl(\Db_M(M)\to\Db_M(U)\bigr),
\]
associating to $U\in\operatorname{Op}_M^\sub$ the subspace of those Schwartz’s distributions on $U$ which extend to the full space. In other words, those distributions that satisfy a tempered growth condition at the boundary $\partial U$.
Note that such a condition on the growth is lost by usual sheafification: $\alpha\Db_M^\temp\simeq\Db_M$.

\begin{remark}
The subanalytic sheaf $\Db_M^\temp$ is a $\beta\D_M$-module. It is not a $\iota\,\D_M$-module, since the analytic coefficients of a differential operator defined on some $U\in\operatorname{Op}_M^\sub$ do not have any growth control at the boundary $\partial U$.
\end{remark}

Here is a real analytic analogue of a fundamental example from \cite[\S 7]{KS03}.

\begin{example}\label{ex:1/x}
On $M=\R_x$, let $P=x^2\partial_x+1$ be the operator associated with the connection $\she_{\R\setminus \{0\}}^{1/x}$, irregular at $x=0$. One has\footnote{To check the isomorphism, note that $e^{1/x}\in \Db_\R^\temp(U)$ if and only if $0\notin \overline{U\cap\{x>0\}}$.}
\[
\ker(\xymatrix{\Db_\R^\temp \ar[r]^-P & \Db_\R^\temp}) \simeq \iota\,\C_{\{x\leq 0\}}\dsum \indlim[\varepsilon\to 0+]\iota\,\C_{\{x> \varepsilon\}}.
\]
This is not a usual sheaf, since one has\footnote{To check the non-isomorphism, compare sections on $U=(0,1)$ using \eqref{eq:ilim}.}
\begin{equation}\label{eq:icolim}
\indlim[\varepsilon\to 0+]\iota\,\C_{\{x> \varepsilon\}} \not\simeq
\iota\mathop{\operatorname{colim}}\limits_{\varepsilon\to 0+}\C_{\{x> \varepsilon\}}
= \iota\,\C_{\{x>0\}}.
\end{equation}
Thus, tempered solutions detect the irregularity of $\she_{\R\setminus \{0\}}^{1/x}$ at $x=0$.
\end{example}

Despite this last example, subanalytic sheaves have some drawbacks in our quest for an irregular Riemann-Hilbert correspondence:

\medskip\noindent \textbf{Caveat II.}\quad
Tempered solutions do not distinguish the connections $\she_{\R\setminus \{0\}}^{1/x}\not\simeq \she_{\R\setminus \{0\}}^{2/x}$, since their associated quoted colimits as in \eqref{eq:icolim} are cofinal.

\medskip\noindent \textbf{Caveat III.}\quad
On $M=\R_t$, let $Q=\partial_t-1$ be the operator associated with the connection $\she_\R^t$, irregular at $t=\infty$. One has
\begin{equation}\label{eq:et}
\ker(\xymatrix{\Db_M^\temp \ar[r]^-Q & \Db_M^\temp}) \simeq \indlim[c\to +\infty]\iota\,\C_{\{t< c\}}\simeq\iota\,\C_M
\quad\in\Mod(\C_\R^\sub).
\end{equation}
Thus, tempered solutions do not detect the irregularity of $\she_\R^t$ at $t=\infty$.

\section{Bordered spaces}
We introduced the notion of bordered spaces in \cite{DK16}, in order to deal with problems like that in Caveat III.

\medskip
A subanalytic bordered space is a pair $M_\infty=(M,\check M)$, where $M$ is an open subanalytic subset of a real analytic manifold $\check M$.
A morphism $r\colon M_\infty \to N_\infty$ of bordered spaces is a morphism $r\colon M\to N$ of real analytic manifolds, such that the closure of its graph
$\overline{\Gamma_r}\subset \check M\times \check N$ is subanalytic, and the projection
$\overline{\Gamma_r}\to \check M$ is proper.
Usual spaces embed fully faithfully in bordered spaces by $M\mapsto (M,M)$.

For example, the identity $\id_M$ and the embedding $M\hookrightarrow\check M$ underlie morphisms
$M \to M_\infty \to \check M$. However, $\id_M$ does not underlie a morphism $M_\infty \to M$, in general.

\medskip
Here is the definition of subanalytic sheaves on bordered spaces:
\begin{itemize}
\item Let $\operatorname{Op}_{M_\infty}^\sub\subset \operatorname{Op}_M^\sub$ be the full subcategory of open subsets which are \emph{relatively compact and subanalytic in $\check M$}.
The category $\Mod(\C_{M_\infty}^\sub)$ of subanalytic sheaves on $M_\infty$ is the full subcategory of presheaves on $\operatorname{Op}_{M_\infty}^\sub$ satisfying the usual gluing conditions on finite covers.
\end{itemize}

As for subanalytic sheaves on usual spaces, there are natural functors
\[
\xymatrix@C=4em{
\Mod(\C_M)
\ar@<.5ex>@{^(->}[r]^-\iota
\ar@<-1.5ex>[r]_-\beta
& \Mod(\C_{M_\infty}^\sub) \ar@<.5ex>[l]|-\alpha .
}
\]
Subanalytic sheaves on bordered spaces have analogous operations as in the case of usual spaces, and we keep the same notations. Formula \eqref{eq:ilim} holds also for bordered spaces:
\begin{equation}\label{eq:ilim2}
\sect(U;\indlim[j]\iota\,F_j)\simeq\mathop{\operatorname{colim}}\limits_j\sect(U;F_j),
\qquad\text{for } U\in\operatorname{Op}_{M_\infty}^\sub.
\end{equation}
For $r\colon M_\infty\to N_\infty$ a morphism, $r_{!!}$ is defined by
\[
\sect(V;r_{!!} F)
= \mathop{\operatorname{colim}}\limits_U\Hom[\Mod(\C_{N_\infty}^\sub)](\iota\,\C_{\opb r V},F\tens\iota\,\C_U),
\qquad\text{for } V\in\operatorname{Op}_{N_\infty}^\sub,
\]
where $U$ ranges over the open subsets in $\operatorname{Op}_{M_\infty}^\sub$ such that $r|_{\opb r V\cap\overline U}$ is proper.

A basic example of bordered space is
\[
\R_\infty \defeq (\R,\mathsf{P}^1), \qquad \mathsf{P}^1=\R\cup\{\infty\}.
\]
Subanalytic sheaves on bordered spaces allow us to strike-out Caveat III:

\medskip\noindent \textbf{\sout{Caveat III}.}\quad
Since\footnote{To check the non-isomorphism, compare sections on an open right half line $U$ using \eqref{eq:ilim2}.} $\indlim[c\to +\infty]\iota\,\C_{\{t< c\}}\not\simeq\iota\,\C_\R$ in $\Mod(\C_{\R_{\infty}}^\sub)$, tempered solutions in $\R_\infty$ do detect the irregularity at $t=\infty$ of $\she_\R^t$.

\section{Enhanced subanalytic sheaves}\label{se:enhsub}

Following \cite{DK16}, let us merge the notions recalled in the previous three sections.

\medskip
Let $M_\infty=(M,\check M)$ be a subanalytic bordered space.
The derived category of enhanced subanalytic sheaves on $M_\infty$ is defined by
\[
\BDC(\C^\enhsub_{M_\infty}) \defeq \{K\in\BDC(\C^\sub_{M_\infty\times\R_\infty})\colon K\ctens E_M^0\isoto K\},
\]
where $\ctens$ is defined using $\reeim\mu$ instead of $\reim\mu$. Consider also the corresponding abelian category (the heart for the natural $t$-structure)
\[
\Mod(\C^\enhsub_{M_\infty}) \defeq \{K\in\BDC(\C^\enhsub_{M_\infty})\colon H^j K\simeq 0 \text{ for }j\neq 0\}.
\]
External operations, relative to a morphism $r\colon M_\infty\to N_\infty$, are defined using external operations for subanalytic sheaves, relative to $\tilde r\defeq r\times\id_{\R_\infty}$. For example, the proper direct image of $K$ is $\reeim{\tilde r}K$.

With notations as in \eqref{eq:EU}, assuming that $\varphi$ underlies a morphism $\varphi\colon U_\infty=(U,\check M)\to \R_\infty$, consider the \emph{exponential} enhanced subanalytic sheaf
\[
{\mathbb E}_{U}^\varphi
\defeq \indlim[c\to +\infty] \iota\, E_U^{\varphi-c}
= \indlim[c\to +\infty] \iota\, \C_{\{t+\varphi(x)\geq c\}}
\quad\in\Mod(\C_{M_\infty}^\enhsub).
\]
Note that ${\mathbb E}_M^0$ is isomorphic to the shifted projection in $\Mod(\C^\enhsub_{M_\infty})$ of the quoted colimit in \sout{Caveat III}:
\begin{equation}\label{eq:c>leq0}
{\mathbb E}_M^0 = \indlim[c\to +\infty] \iota\, \C_{\{t\geq c\}} \simeq \indlim[c\to +\infty] \iota\, \C_{\{t< c\}}[1] \ctens E_M^0.
\end{equation}

Consider the fully faithful embeddings, compatible with operations:
\[
\xymatrix{
\BDC(\C_M) \ar@{^(->}[r]^-\iota
& \BDC(\C^\sub_{M_\infty}) \ar@{^(->}[r]^-e
& \BDC(\C^\enhsub_{M_\infty}),
}
\]
with $e\,G\defeq\opb\pi G\ctens {\mathbb E}_M^0$.
Note that $e\,\iota\, F\simeq (\iota\,\epsilon\, F)\ctens {\mathbb E}_M^0$.

\medskip
Enhanced subanalytic sheaves provide the correct setting to stage the irregular Riemann-Hilbert correspondence.
For example, we can strike-out the two remaining caveats:

\medskip\noindent \textbf{\sout{Caveat II}.}\quad
${\mathbb E}_{\R\setminus\{0\}}^{1/x}\not\simeq {\mathbb E}_{\R\setminus\{0\}}^{2/x}$ as enhanced subanalytic sheaves on $\R$.

\medskip\noindent \textbf{\sout{Caveat I}.}\quad
${\mathbb E}_X^{\Re a}\simeq {\mathbb E}_X^0$ as enhanced subanalytic sheaves on $X$.

\section{Back to the complex domain}

Let $X$ be  a complex manifold, and $X_\R$ the underlying real analytic manifold.
Recall that the sheaf $\O_X$ is quasi-isomorphic to the Dolbeault complex with coefficients in $\Db_{X_\R}$.  Following \cite{KS01}, one sets
\begin{itemize}
\item $\O_X^\temp = \left(\text{Dolbeault complex with coefficients in }\Db_{X_\R}^\temp\right)$,
\item $\sol^\temp(\shm) = \rhom[\beta\D_X](\beta\shm,\O_X^\temp)$.
\end{itemize}
Expressed in this terminology, the left pointing arrow in Theorem~\ref{th:Kas84} is given by
\begin{equation}\label{eq:thom}
\BDC(\C_X)\to\BDC(\D_X),\quad
F\mapsto\rhom(\iota\,F,\O_X^\temp),
\end{equation}
using the ``stacky'' hom for subanalyic sheaves, with values in usual sheaves.

\medskip
Let $\mathbb{P}^1=\C_\tau\cup\{\infty\}$, with $t=\Re\tau$ the affine coordinate of $\mathsf P^1=\R_t\cup\{\infty\}$. Consider the natural composition
$i\colon X\times\R_\infty\to X\times \mathsf{P}^1\to X\times \mathbb{P}^1$.
In \cite{DK16} we introduced the following enhancements\footnote{The shift $[2]$ in the definition of $\O_X^\Temp$ is motivated by a first shift from $\Db^\temp_{X\times\R_\infty}\simeq\epb i\O^\temp_{X\times\mathbb{P}^1}[1]$, and a second shift from \eqref{eq:c>leq0}.}:
\begin{itemize}
\item $\O_X^\Temp = \epb i \rhom[\beta\D_{\mathbb{P}^1}](\beta\she_{\mathbb{P}^1}^\tau,\O_{X\times\mathbb{P}^1}^\temp)[2] \ctens E_X^0$,
\item $\sol^\Temp(\shm) = \rhom[\beta\D_X](\beta\shm,\O_X^\Temp)
\simeq  \epb i \sol^\temp(\shm\etens_\O\she_{\mathbb{P}^1}^\tau)[2]  \ctens E_X^0$.
\end{itemize}
Note that $\O_X^\Temp \simeq \O_X^\Temp \ctens \mathbb{E}^0_X$.
We also introduced the natural enhancement of \eqref{eq:thom}
\[
\BDC(\C_X^\enhsub)\to\BDC(\D_X),\quad
K\mapsto\rhom(K,\O_X^\Temp),
\]
which satisfies $\rhom(e\,\iota\,F,\O_X^\Temp) \simeq \rhom(\iota\,F,\O_X^\temp)$.

\begin{proposition}
For\footnote{As usual, we denote with an index, like in $\BDC_\reghol(\D_X)$, the full triangulated subcategory of objects whose cohomologies satisfy that property. In this case, complexes of $\D_X$-modules with regular holonomic cohomologies.} $\shl\in\BDC_\reghol(\D_X)$, one has $\sol^\Temp(\shl) \simeq e\,\iota\,\sol(\shl)$.
\end{proposition}

This is an enhancement of the classical result $\sol^\temp(\shl)
\simeq\iota\,\sol(\shl)$.

\begin{proposition}\label{pro:solT}
For $f\in\O_X(*D)$, one has
\begin{itemize}
\item[(i)]
$\sol^\Temp(\she_{X\setminus D}^f) \simeq \mathbb{E}_{X\setminus D}^{\Re f}$,
\item[(ii)]
$\rhom(\mathbb{E}_{X\setminus D}^{\Re f},\,\O_X^\Temp) \simeq \she_{X\setminus D}^f$.
\end{itemize}
\end{proposition}

Here: (i) can be reduced to a computation like in Example~\ref{ex:1/x}, followed by a projection like in \eqref{eq:c>leq0}; and (ii) is analogous to \cite[Proposition~8.1]{Dag14}.

\section{The Riemann-Hilbert correspondence}\label{sec:RH}

In the theorem below, the top row expresses our irregular Riemann-Hilbert correspondence from \cite{DK16}. It is compatible with the regular Riemann-Hilbert correspondence from \cite{Kas84}, reproduced in the bottom row\footnote{The index in $\BDC_\Cc(\C_X)$ stands for sheaves which are constructible with respect to a complex analytic stratification.}.

\begin{theorem}\label{th:RHregirr}
Let $X$ be  a complex manifold.
There is a commutative diagram, compatible with Grothendieck's operations,
\[
\xymatrix@C=4em{
\BDC_\hol(\D_X) \ar@{^(->}[r]_-{\sol^\Temp} \ar@{^(->}@/^2em/[rrr]^\empty
& \BDC(\C_X^\enhsub)  \ar[rr]_-{\rhom(\bullet,\,\O_X^\Temp)}
&& \BDC(\D_X) \\
\BDC_\reghol(\D_X) \ar[r]_-{\sol}^-\sim \ar@{^(->}[u]
& \BDC_\Cc(\C_X)  \ar[rr]_-{\rhom(\iota(\bullet),\,\O_X^\temp)}^-\sim   \ar@{^(->}[u]^{e\circ\iota}
&& \BDC_\reghol(\D_X) , \ar@{^(->}[u]
}
\]
where the unnamed arrows are inclusions of full subcategories.
\end{theorem}

In particular, as the arrow $\hookrightarrow$ indicates, we proved that the enhanced tempered solution functor $\sol^\Temp$ is fully faithful. And we also proved that a holonomic $\D$-module $\shm$ can be reconstructed from its enhanced tempered solutions $\sol^\Temp(\shm) \eqdef K$, by $\shm\simeq\rhom(K,\O_X^\Temp)$.

Our proof makes full use of dévissage and functoriality, like it was the case in \cite{Kas84}. There, Kashiwara reduced to a normal form for regular flat meromorphic connections established by Deligne~\cite{Del70}. Here, we reduce to a normal form for flat meromorphic connections established in the deep works of Mochizuki~\cite{Moc11} and Kedlaya~\cite{Ked11}, following a conjecture and preliminary results by Sabbah~\cite{Sab00}. After such a reduction, Proposition~\ref{pro:solT} can be put to good use.

As mentioned in Footnote \ref{foot:disc}, a topological description of the essential image of $\sol^\Temp$ is still missing. That is what one should call $\BDC_\Cc(\C_X^\enhsub)$, by analogy with $\BDC_\Cc(\C_X)$. In \cite{DK16} we did a first step in this direction, by showing that $\sol^\Temp$ takes values in the category $\BDC_\Rc(\C_X^\enhsub)$. This is the enhanced subanalytic analogue of the usual $\BDC_\Rc(\C_X)$, where the index stands for sheaves which are constructible with respect to a subanalytic stratification.

The abelian analogue of the top row in Theorem~\ref{th:RHregirr} is given in Theorem~\ref{th:RHirr}.
Again, a topological description of perverse enhanced subanalytic sheaves is still missing. In \cite{DK19} we did a first step in this direction, by introducing a notion of perverse $t$-structure in $\BDC_\Rc(\C_X^\enhsub)$.

A curve test for the essential image was established by Mochizuki~\cite{Moc22}. It states that a complex $K$ of enhanced subanalytic sheaves belongs to the essential image of $\Mod_\hol(\D_X)$ by $\sol^\Temp$ if and only if the restriction of $K$ to any smooth analytic curve satisfies the same condition. The interest of this criterion lies in the fact that, in the one dimensional case, one has a complete description of the essential image of $\sol^\Temp$. We present such description in the next section, for the case of meromorphic connections.

\section{Stokes local systems}\label{sec:SLS}
Let us upgrade the ``Stokes local systems'' from Example~\ref{ex:Airy}.
This will provide a target category for an irregular analogue of Theorem~\ref{th:Del}, in the one dimensional case.

\medskip
Let $X$ be a smooth complex curve, and $D\subset X$ a discrete set of points. The set of exponential factors at $c\in D$ is
\[
\operatorname{Expf_c} \defeq \mathop{\operatorname{colim}}\limits_{d\in\Z_{>0}} z_c^{-1/d}\,\C[z_c^{-1/d}],
\]
where $z_c$ is a local coordinate at $c$ with $z_c(c)=0$, and $\Z_{>0}$ is partially ordered by divisibility.

Let $K\in\Mod(\C_{(X\setminus D)_\infty}^\enhsub)$,
with $(X\setminus D)_\infty \defeq (X\setminus D, X)$.
One says that:
\begin{itemize}
\item $K$ is a local system on $X\setminus D$, if $\opb{\tilde\imath} K\simeq e\,\iota\,L$. Here, $L$ is a local system on $X\setminus D$, and $i\colon X\setminus D\to (X\setminus D)_\infty$ is the natural morphism;
\item $K$ is of Stokes type at $c$, if there are $f_1,\dots,f_m\in\operatorname{Expf_c}$, and a finite open sectorial cover $\{V_\theta\}_{\theta\in I}$ of a small neighborhood of $c$ in $X\setminus D$, such that for any $\theta\in I$
\[
\opb{\tilde\jmath_\theta}K\simeq \DSum_{j=1}^m{\mathbb E}_{V_\theta}^{\Re f_j}\quad\in\Mod(\C^\enhsub_{V_{\theta,\infty}}).
\]
Here, $j_\theta\colon V_{\theta,\infty}=(V_\theta,X)\to (X\setminus D)_\infty$ is the natural morphism, and we fixed a determination of each $f_j$ on each $V_\theta$.
\end{itemize}

The category of Stokes local systems on $(X\setminus D)_\infty$ is the full abelian subcategory of $\Mod(\C^\enhsub_{(X\setminus D)_\infty})$ whose objects are both local systems on $X\setminus D$, and of Stokes type at any $c\in D$.

\begin{theorem}\label{th:SLS}
There is an equivalence of abelian categories
\[
\xymatrix@C=6em{
\left\{
\text{$\O_X(*D)$-connections}
\right\}
\ar[r]^-{\sol^\Temp|_{(X\setminus D)_\infty}}_-\sim
&
\left\{
\text{Stokes local systems on $(X\setminus D)_\infty$}
\right\}^\op
}.
\]
\end{theorem}

In the late 1980's Deligne and Malgrange (see \cite{DMR07}) provided an alternative equivalence, whose target is the category of Stokes filtered local systems\footnote{Such local systems live on the real oriented blow-up of $X$ along $D$.}. The essential surjectivity of $\sol^\Temp|_{(X\setminus D)_\infty}$ in Theorem~\ref{th:SLS} is obtained by naturally associating a Stokes filtered local system to a Stokes local system.

\address{Dipartimento di Matematica, Universit{\`a} di Padova\\
via Trieste 63, 35121 Padova, Italy\\
\email{dagnolo@math.unipd.it}}

\address{Kyoto University Institute for Advanced study\\
Research Institute for Mathematical Sciences, Kyoto University\\
Kyoto 606-8502, Japan\\
\email{masaki@kurims.kyoto-u.ac.jp}}

\end{document}